\documentclass[11pt,reqno]{article}
\usepackage{amsmath,mathrsfs,graphicx,amssymb,amsfonts,color}

\font\bg=cmbx10 scaled\magstep1

\font\small=cmr8

\newtheorem{newlemma}{{\bf Lemma}}

\newtheorem{newteorem}{{\bf Theorem}}

\newenvironment{teorem}{\begin{newteorem}{\hspace{-0.5
em}{\bf.}}}{\end{newteorem}}

\newtheorem{newkorolari}{{\bf Corollary}}

\newtheorem{newdefine}{{\bf Definition}}

\newtheorem{newquestion}{{\bf Question}}

\newtheorem{newkonjek}{{\bf Conjecture}}

\newtheorem{newexample}{{\bf Example}}

\begin{document}
\tolerance=10000
\baselineskip18truept
\newbox\thebox
\global\setbox\thebox=\vbox to 0.2truecm{\hsize 0.15truecm\noindent\hfill}
\def\boxit#1{\vbox{\hrule\hbox{\vrule\kern0pt
\vbox{\kern0pt#1\kern0pt}\kern0pt\vrule}\hrule}}
\def\qed{\lower0.1cm\hbox{\noindent \boxit{\copy\thebox}}\bigskip}
\def\ss{\smallskip}
\def\ms{\medskip}
\def\bs{\bigskip}
\def\c{\centerline}
\def\nt{\noindent}
\def\ul{\underline}
\def\ol{\overline}
\def\lc{\lceil}
\def\rc{\rceil}
\def\lf{\lfloor}
\def\rf{\rfloor}
\def\ov{\over}
\def\t{\tau}
\def\th{\theta}
\def\k{\kappa}
\def\l{\lambda}
\def\L{\Lambda}
\def\g{\gamma}
\def\d{\delta}
\def\D{\Delta}
\def\e{\epsilon}
\def\lg{\langle}
\def\rg{\rangle}
\def\p{\prime}
\def\sg{\sigma}
\def\ch{\choose}

\newcommand{\ben}{\begin{enumerate}}
\newcommand{\een}{\end{enumerate}}
\newcommand{\bit}{\begin{itemize}}
\newcommand{\eit}{\end{itemize}}
\newcommand{\bea}{\begin{eqnarray*}}
\newcommand{\eea}{\end{eqnarray*}}
\newcommand{\bear}{\begin{eqnarray}}
\newcommand{\eear}{\end{eqnarray}}

\centerline{\Large \bf  Domination polynomial of generalized friendship   }
\vspace{.3cm}

\centerline {\Large \bf and generalized book  graphs }

\bigskip

\bs

\baselineskip12truept
\centerline{S. Jahari$^{}${}\footnote{\baselineskip12truept\it\small
Corresponding author. E-mail: s.jahari@gmail.com} and S. Alikhani }
\baselineskip20truept
\centerline{\it Department of Mathematics, Yazd University}
\vskip-8truept
\centerline{\it  89195-741, Yazd, Iran}

\vskip-0.2truecm
\nt\rule{16cm}{0.1mm}

\nt{\bg ABSTRACT}
\medskip

\baselineskip14truept

\nt{ Let $G$ be a simple graph of order $n$.
The domination polynomial of $G$ is the polynomial
$D(G, x)=\sum_{i=\gamma(G)}^{n} d(G,i) x^{i}$,
where $d(G,i)$ is the number of dominating sets of $G$ of size $i$ and
$\gamma(G)$ is the domination number of $G$.
 Let $n$ be any positive integer and  $F_n$ be the {\em Friendship graph} with
$2n + 1$ vertices and $3n$ edges, formed by the join of $K_{1}$ with $nK_{2}$.
We study the domination polynomials of generalized friendship  graphs.
We also consider the {\em $n$-book graphs } $B_n$, formed by joining $n$ copies of the cycle graph $C_4$ with a common edge and study the domination polynomials of some generalized book  graphs. In particular we examine the domination roots of these families, and find the limiting curve for the roots.}


\ms

\nt{\bf Mathematics Subject Classification:} {\small 05C60.}

\nt\rule{16cm}{0.1mm}

\baselineskip20truept

\section{Introduction}

\nt Let $G=(V,E)$ be a simple graph.
For any vertex $v\in V(G)$, the {\it open neighborhood} of $v$ is the
set $N(v)=\{u \in V (G) | \{u, v\}\in E(G)\}$ and the {\it closed neighborhood} of $v$
is the set $N[v]=N(v)\cup \{v\}$. For a set $S\subseteq V(G)$, the open
neighborhood of $S$ is $N(S)=\bigcup_{v\in S} N(v)$ and the closed neighborhood of $S$
is $N[S]=N(S)\cup S$.
A set $S\subseteq V(G)$ is a {\it dominating set} if $N[S]=V$ or equivalently,
every vertex in $V(G)\backslash S$ is adjacent to at least one vertex in $S$.
The {\it domination number} $\gamma(G)$ is the minimum cardinality of a dominating set in $G$.
For a detailed treatment of these parameters, the reader is referred to~\cite{domination}.
The {\it $i$-subset} of $V(G)$ is a subset of $V(G)$ of size $i$.
Let ${\cal D}(G,i)$ be the family of dominating sets of a graph $G$ with cardinality $i$ and
let $d(G,i)=|{\cal D}(G,i)|$.
The {\it domination polynomial} $D(G,x)$ of $G$ is defined as
$D(G,x)=\sum_{ i=\gamma(G)}^{|V(G)|} d(G,i) x^{i}$,
where $\gamma(G)$ is the domination number of $G$ (see \cite{euro,saeid1}). A 
 root of $D(G, x)$ is called a domination root of $G$. The set of distinct roots of $D(G, x)$ is denoted by $Z(D(G, x))$.

\nt Calculating the domination polynomial of a graph $G$ is difficult in general, as the smallest power of a non-zero term is the domination number $\gamma (G)$ of the graph, and determining whether $\gamma (G) \leq k$ is known to be NP-complete \cite{garey}. But for certain classes of graphs, we can find a closed form expression for the domination polynomial. In \cite{jason} the domination polynomial and the domination roots of friendship graphs has been studied. In this paper we would like to obtain some further results of this kind. We consider generalized friendship graph (or flower graphs), and generalized  book graphs and calculate their domination polynomials. Also we  explore  the nature and location of their roots.

\section{Domination polynomial of generalized friendship graphs}

\begin{figure}[ht]
\hspace{3.6cm}
\includegraphics[width=8.5cm,height=2cm]{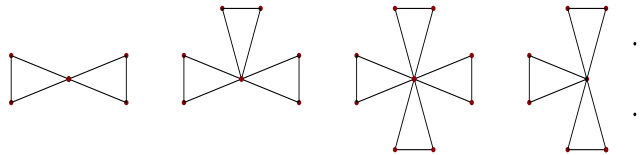}
\caption{\label{dutch} Friendship graphs $F_{2}, F_{3}, F_{4}$ and $F_{n}$, respectively. }
\end{figure}

\nt Let consider the graphs $F_{n}$ obtained by selecting one vertex in each of $n$ triangles and identifying them (Figure \ref{dutch}). Some call them Dutch-Windmill graphs \cite{htt} and friendship graphs.

\begin{figure}[ht]
\hspace{3.6cm}
\includegraphics[width=8.5cm,height=2cm]{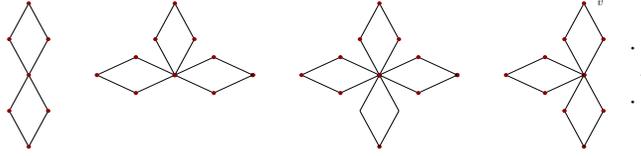}
\caption{\label{figtetragons} The flowers  $F_{4,2},~ F_{4,3}, ~F_{4,4}$ and $F_{4,n}$, respectively. }
\end{figure}

\nt The generalized friendship graph $F_{q,n}$ is a collection of $n$ cycles (all of order $q$), meeting at a common vertex (see Figure \ref{figtetragons}). The generalized friendship graph may also be referred
to as a flower \rm\cite{schi}.

\nt In this section we compute the domination polynomial of the flowers  $F_{4,n}$. We need some preliminaries.

\nt The vertex contraction $G/u$ of a graph $G$ by a vertex $u$ is the operation under
which all vertices in $N(u)$ are joined to each other and then $u$ is deleted (see\cite{Wal}).

\nt The following theorem is  useful for finding the recurrence relations for the  domination polynomials  of  graphs.

\begin{teorem}\label{theorem1}{\rm \cite{saeid2,Kot}}
Let $G$ be a graph. For any vertex $u$ in $G$ we have
\[
D(G, x) = xD(G/u, x) + D(G - u, x) + xD(G - N[u], x) - (1 + x)p_u(G, x),
\]
where $p_u(G, x)$ is the polynomial counting the dominating sets of $G - u$ which do not contain any
vertex of $N(u)$ in $G$.
\end{teorem}

\nt The following theorem gives formula for the domination polynomial of $F_{n}$.
\begin{teorem} \label{theorem4}{\rm \cite{jason}}
For every $n\in \mathbb{N}$,
$$D(F_{n},x) = (2x + x^2)^n +x(1 + x)^{2n}.$$
\end{teorem}

\nt Domination polynomial satisfies a recurrence relation for arbitrary graphs which is based on the edge and vertex
elimination operations. The recurrence uses composite operations, e.g. $G - e /u$, which
stands for $(G - e) /u$.

\begin{teorem}\label{theorem3}{\rm \cite{ Kot}}
Let $G$ be a graph. For every edge $e =\{u, v\}\in E$,
\begin{eqnarray*}
D(G, x)&=&D(G - e, x) +\frac{x}{x-1} \Big[ D(G - e/u, x) + D(G - e/v, x)\\
&-&D(G/u, x) - D(G/v, x) - D(G - N[u], x) - D(G - N[v], x)\\
&+&D(G - e - N[u], x) + D(G - e - N[v], x)\Big].
\end{eqnarray*}
\end{teorem}


\nt The following theorem gives recurrence relation for the domination polynomial of $F_{4,n}$.

\begin{teorem}\label{theorem5}
\nt For every $n \geq 2$,
\begin{eqnarray*}
D(F_{4,n},x)&=&((1+x)^3 +x)D(F_{4,n-1},x)- (1 +3x)(x+3x^2+x^3)^{n-1}\\
 && +(1 + x)^3x^{n-1}- (x^2+x)(x^3+3x^2+3x)^{n-1},
\end{eqnarray*}
where $D(F_{4,1},x) = x^4+4x^3+6x^2.$
\end{teorem}
\nt{\bf Proof.} An elementary observation is that if $G_1$ and $G_2$ are graphs of orders $n_1$ and $n_2$,
respectively, then
\[ D(G_1 \cup G_2,x) = D(G_1, x) D(G_2, x).\]
\nt Consider graph  $F_{4,n}$ and a vertex $v$ in  Figure \ref{figtetragons}. By Theorem ~\ref{theorem1} we have:
\begin{eqnarray}\label{eq1}
D(F_{4,n}, x)&=& x D(F_{4,n}/v, x) + D(F_{4,n} - v, x) + x D(F_{4,n} - N[v], x) - (1 + x)p_v(F_{4,n}, x)\nonumber\\
&=& x D(F_{4,n}/v, x) + D(F_{4,n} - v, x) + xD(F_{4,n-1} ,x)\nonumber\\
&& - (1 + x)xD(\cup_{i=1}^{n-1}K_3 ,x),
\end{eqnarray}
where $D(K_3,x)=x^3+3x^2+3x$.

\begin{figure}[ht]
\hspace{2.3cm}
\includegraphics[width=9.5cm,height=2.9cm]{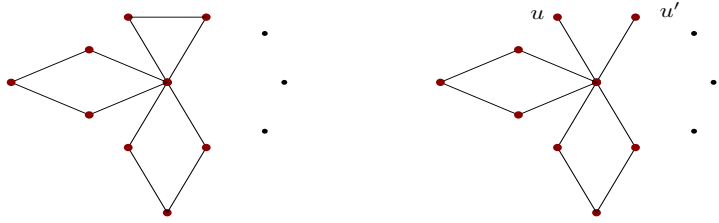}
\caption{ \label{figure1} Graphs $F_{4,n}/v$ ~and ~ $F_{4,n} - v$,~ respectively.}
\end{figure}

\nt Now  we use Theorem~\ref{theorem3} to obtain the domination polynomial of the graph $F_{4,n}/v=G$
 (see Figure~\ref{figure1}). We have

\begin{eqnarray}\label{eq2}
D(F_{4,n}/v, x)&=&D(F_{4,n} - v, x) +\frac{x}{x-1} \Big[ - D(\cup_{i=1}^{n-1}P_3 ,x) - D(\cup_{i=1}^{n-1}P_3 ,x)\nonumber\\
&+&xD(\cup_{i=1}^{n-1}P_3 ,x) + xD(\cup_{i=1}^{n-1}P_3 ,x)\Big]\nonumber\\
&=&D(F_{4,n} - v, x) +2xD(P_3 ,x)^{n-1},
\end{eqnarray}

\nt where $ D(G - e/u, x) \simeq D(G - e/u', x) \simeq D(G/u, x) \simeq D(G/u', x) $ and $(F_{4,n}/v) - e\simeq F_{4,n} - v$ in Figure~\ref{figure1}.

\nt Now  we use Theorem~\ref{theorem1} to obtain the domination polynomial of the graph $F_{4,n} - v=G$
 (see Figure~\ref{figure1}). We have

\begin{eqnarray}\label{eq3}
D(F_{4,n} - v, x)&=& x D((F_{4,n} - v)/u, x) + D(F_{4,n} - v- u, x) + x D(F_{4,n} - v- N[u], x)  \nonumber\\
&&- (1 + x)p_u(F_{4,n} - v, x)\nonumber\\
&=&x D((F_{4,n} - v)/u, x) + D(F_{4,n} - v- u, x) + x(xD(\cup_{i=1}^{n-1}P_3 ,x)) \nonumber\\
&&- (1 + x)xD(\cup_{i=1}^{n-1}P_3 ,x)\nonumber\\
&=&(1+x) D((F_{4,n} - v)/u, x) - xD(P_3 ,x)^{n-1}.
\end{eqnarray}

\begin{figure}[ht]
\hspace{5.5cm}
\includegraphics[width=3.6cm,height=2.4cm]{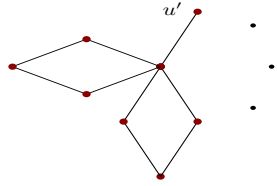}
\caption{ \label{figure2} The graph $(F_{4,n} - v)/u$.}
\end{figure}

\nt Note that we used  $ D((F_{4,n} - v) /u, x) \simeq D(F_{4,n} - v - u, x) $  (see Figure~\ref{figure2}).

\nt Use Theorem~\ref{theorem1} to obtain the domination polynomial of the graph $(F_{4,n} - v)/u=G$
 (see Figure~\ref{figure2}). We have

\begin{eqnarray}\label{eq4}
D(F_{4,n} - v)/u, x)&=& x D((F_{4,n} - v)/u', x) + D(F_{4,n} - v - u', x) + x D(F_{4,n} - v- N[u'], x) \nonumber\\
&&- (1 + x)p_{u'}(F_{4,n} - v, x)\nonumber\\
&=&x D(F_{4,n-1} ,x) + D(F_{4,n-1} ,x) + xD(\cup_{i=1}^{n-1}P_3 ,x)\nonumber\\
&&- (1 + x)(D(\cup_{i=1}^{n-1}P_3 ,x)-x^{n-1})\nonumber\\
&=&(1+x) D(F_{4,n-1} ,x) - D(P_3 ,x)^{n-1} +(1 + x) (x^{n-1}) ,
\end{eqnarray}

\nt Consequently, by equations \ref{eq1}, \ref{eq2}, \ref{eq3} and \ref{eq4} have:
\begin{eqnarray*}
D(F_{4,n},x)&=&(1+4x+3x^2+x^3)D(F_{4,n-1},x)- (1 +3x)(x+3x^2+x^3)^{n-1}\\
 && +(1 + x)^3x^{n-1}- (x^2+x)(x^3+3x^2+3x)^{n-1}.\quad\qed
\end{eqnarray*}

\begin{figure}
\hspace{5cm}
\includegraphics[width=7cm]{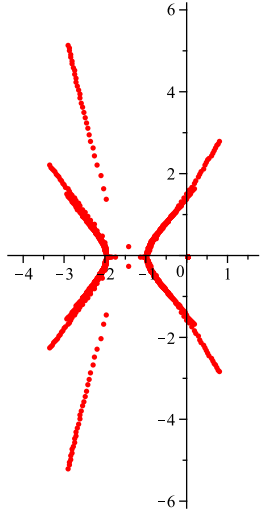}
\caption{\label{figure2'} Domination roots of graphs $F_{4,n}$, for $1 \leq n \leq 30$.}
\end{figure}

\nt The domination roots of $F_{4,n}$ exhibit a number of interesting properties (see Figure~\ref{figure2'}).

\nt If we can find an explicit formula for the domination polynomial of a graph, there are still interesting, difficult problems concerning the roots. For every odd natural number $n$, no nonzero real number is a domination
root of $F_{n}$ \cite{jason}. Also we think that for $n$ even,  $F_{n}$ have exactly three real roots. Using Maple, we think that these are  true for $F_{4,n}$. Therefore,  we pose the following:

\begin{newquestion}
For even $n\geq 4$, does $F_{n}$ have exactly three real roots?
\end{newquestion}

\begin{newkonjek}
For every odd natural number $n$, no nonzero real number is a domination
root of $F_{4,n}$.
\end{newkonjek}

\begin{newquestion}
What is a good upper bound on the modulus of the roots of  $F_{4,n}$?
\end{newquestion}

\nt It is natural to ask about the complex domination roots of $F_{4,n}$.  The plot in Figure~\ref{figure2'} suggests that the roots tend to lie on a curve.
\begin{newkonjek}
 The limit of domination roots of $F_{4,n}$ is hyperbola.
 \end{newkonjek}

\nt In \cite{brown} a family of graphs was produced with roots just barely in the right-half plane (showing that not all domination polynomials are stable), but Figure~\ref{figure2'}  provides an explicit family (namely the $F_{4,n}$) whose domination roots have unbounded positive real part.
\ms

\section{Domination polynomial of generalized book graphs}

\begin{figure}[ht]
\hspace{2.3cm}
\includegraphics[width=9.5cm,height=2.9cm]{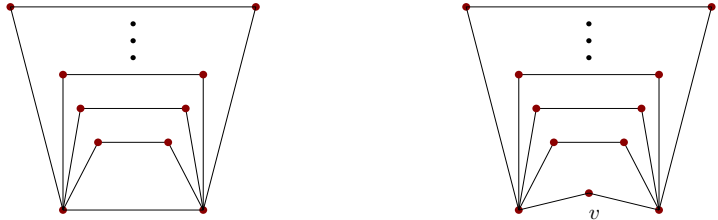}
\caption{ \label{figure3} Graphs $B_n$ ~and ~ $B_{n,5} $,~ respectively.}
\end{figure}

\nt A book graph $B_n$, is defined as follows
$V(B_n)=\{u_1,u_2\}\cup \{v_i, w_i : 1\leq i\leq n\}$ and $E(B_n)=\{u_1u_2\}\cup \{u_1v_i,~ u_2w_i,~v_iw_i : 1\leq i\leq n\}$.
We consider  the generalized book graph $B_{n,m}$ with vertex and edge sets by
$V(B_{n,m})=\{u_i: 1\leq i \leq m-2\}\cup \{v_i, w_i : 1\leq i\leq n\}$ and $E(B_{n,m})=\{u_iu_{i+1}:   1\leq i \leq m-3\}\cup \{u_iw_j : 1\leq j\leq n,~i=m-2\}\cup \{ u_1v_i : 1\leq i\leq n\}\cup \{v_iw_i : 1\leq i\leq n\}$ (see Figure \ref{figure3}).

\nt The following theorem gives formula for the domination polynomial of $B_n$.

\begin{teorem}\label{theorem6}{\rm \cite{jason}}
\nt For every $n \in \mathbb{N}$,
\[ D(B_n,x)=(x^2+2 x)^n(2x+1) + x^2(x+1)^{2n}- 2x^n.\]
\end{teorem}

\begin{figure}[ht]
\hspace{4cm}
\includegraphics[width=7cm]{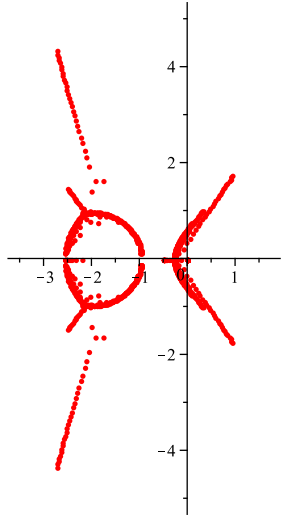}
\caption{\label{bookroots} Domination roots of graphs $B_{n}$, for $1 \leq n \leq 30$.}
\end{figure}

\nt Figure~\ref{bookroots} shows the domination roots of book graphs $B_{n}$ for $n \leq 30$.

\nt In this section we compute domination polynomial of the book graphs  $B_{n,5}$. We need some preliminaries.

\nt We begin  with a graph operation.  For two graphs $G = (V,E)$ and $H=(W,F)$, the corona $G\circ H$ is the graph arising from the
disjoint union of $G$ with $| V |$ copies of $H$, by adding edges between
the $i$th vertex of $G$ and all vertices of $i$th copy of $H$ \cite{Fruc}. It is easy to see that the corona operation of two graphs does not have the commutative property.
The following theorem  gives  the domination polynomial of graphs of
the form $H\circ K_1$, which is needed to obtain our result.

\begin{teorem}\label{theorem8}{\rm \cite{euro}}
Let $G $  be a graph. Then $D(G,x) = (x^2 + 2x)^n$ if and only if $G=H\circ K_1$  for some graph  $H$ of order $n$.
\end{teorem}

\nt Given any two graphs $G$ and $H$ we define the {\it Cartesian product}, denoted $G\Box H$, to be the graph with vertex set $V(G)\times V(H)$ and edges between two vertices $(u_1,v_1)$ and  $(u_2,v_2)$ if and only if either $u_1=u_2$ and $v_1v_2 \in E(H)$ or $u_1u_2 \in E(G)$ and $v_1=v_2$. This product is well known to be commutative.

 \begin{teorem}\label{theorem9}{\rm \cite{Kotek}}
The domination polynomial for $K_n\Box K_2$ is
\[
D(K_n\Box K_2,x)=((1+x)^n-1)^2+2x^n.
\]
\end{teorem}

\nt The following theorem gives formula for the domination polynomial of $B_{n,5}$.

\begin{teorem}\label{theorem7}
\nt For every $n \in \mathbb{N}$,
\begin{eqnarray*}
D(B_{n,5},x)&=x^2(x+1)^{2n+1} - 2x^{n+1} + (x^2+2x)^{n}(2x^2+3x).
\end{eqnarray*}
\end{teorem}
\nt{\bf Proof.}
\nt Consider graph  $B_{n,5}$ in Figure \ref{figure3}. By Theorems ~\ref{theorem1} and  \ref{theorem6} we have:
\begin{eqnarray}\label{eq5}
D(B_{n,5}, x)&=& x D(B_{n,5}/v, x) + D(B_{n,5} - v, x) + x D(B_{n,5} - N[v], x)\nonumber\\
&& - (1 + x)p_v(B_{n,5}, x)\nonumber\\
&=& x D(B_n, x) + D(B_{n,5} - v, x) + x(D(\cup_{i=1}^{n}K_2,x))\nonumber\\
&& - (1 + x)[(x^2+2x)^n - 2x^n]\nonumber\\
&=&x[(x^2+2 x)^n(2x+1) + x^2(x+1)^{2n}- 2x^n] +D(B_{n,5} - v, x) \nonumber\\
&& +  x(x^2+2x)^{n} - (1 + x)[(x^2+2x)^n - 2x^n]\nonumber\\
&=& x^3(x+1)^{2n} + 2x^{n} +D(B_{n,5} - v, x) +   (x^2+2x)^{n}(2x^2+x-1).
\end{eqnarray}

\begin{figure}[ht]
\hspace{5.5cm}
\includegraphics[width=5cm,height=2.9cm]{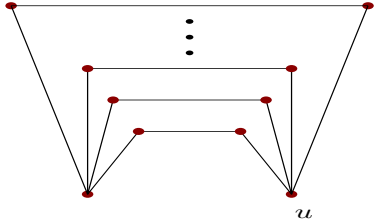}
\caption{ \label{figure4} The graph $B_{n,5} - v$.}
\end{figure}

\nt Now  we use Theorem~\ref{theorem1} to obtain the domination polynomial of the graph $B_{n,5} - v=G$
 (see Figure~\ref{figure4}).  We have

\begin{eqnarray}\label{eq6}
D(B_{n,5} - v, x)&=& x D(G/u, x) + D(G - u, x) + x D(G - N[u], x) - (1 + x)p_u(G, x)\nonumber\\
&=&x D(G/u, x) + D(G - u, x) + xD(K_{1,n} ,x) - (1 + x)(x^n(1+x))\nonumber\\
&=&x D(G/u, x) + D(G - u, x) +x(x^n+x(1+x)^n) - x^n(1+x)^2,
\end{eqnarray}

\nt where $ D(K_{1,n} ,x) = x^n+x(1+x)^n $.

\begin{figure}[ht]
\hspace{3.5cm}
\includegraphics[width=9.4cm,height=2.9cm]{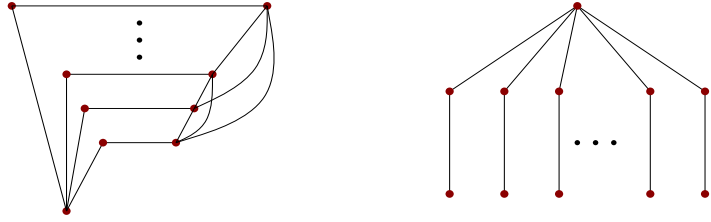}
\caption{ \label{figure5} Graphs $(B_{n,5} - v)/u$ and $B_{n,5} - v - u$,~ respectively.}
\end{figure}

\nt Use Theorems~\ref{theorem1} and \ref{theorem8} to obtain the domination polynomial of the graph $B_{n,5} - v - u=G$
 (see Figure~\ref{figure5}). We have

\begin{eqnarray}\label{eq7}
D(B_{n,5} - v - u, x)&=& x D(G/w, x) + D(G - w, x) + x D(G- N[w], x)\nonumber\\
&& - (1 + x)p_{w}(G, x)\nonumber\\
&=&x D(K_n\circ K_1 ,x) + D(\cup_{i=1}^{n}P_2 ,x)+x(x^n)- (1 + x)x^n\nonumber\\
&=& (2x+x^2)^n(x+1) - x^n.
\end{eqnarray}

\nt Use Theorems~\ref{theorem1} and \ref{theorem8} and \ref{theorem9} to obtain the domination polynomial of the graph $(B_{n,5} - v)/u=G$
 (see Figure~\ref{figure5}). We have

\begin{eqnarray}\label{eq8}
D((B_{n,5} - v)/u, x)&=& x D(G/u', x) + D(G - u', x) + x D(G- N[u'], x) - (1 + x)p_{u'}(G, x)\nonumber\\
&=&x D(K_n\Box K_2 ,x) + D(K_n\circ K_1 ,x) + xD(K_n ,x) - (1 + x)x^n\nonumber\\
&=&x((1+x)^n-1)^2 +2x^n) + (2x+x^2)^n + x((1+x)^n-1) - (1+x)x^n\nonumber\\
&=&  x((1+x)^n-1)(1+x)^n+ (2x+x^2)^n + x^n (x-1).
\end{eqnarray}

\nt By equations \ref{eq6}, \ref{eq7} and \ref{eq8} have:
\begin{eqnarray}\label{eq9}
D(B_{n,5} - v,x)&=&(x^2+2 x)^n(2x+1) + x^2(x+1)^{2n} - 2x^n(1 + x).
\end{eqnarray}

\nt Consequently, by equations \ref{eq5} and \ref{eq9} have:
\begin{eqnarray*}
D(B_{n,5},x)&=& x^3(x+1)^{2n} + 2x^{n} +   (x^2+2x)^{n}(2x^2+x-1) \\
&& + (x^2+2 x)^n(2x+1) + x^2(x+1)^{2n} -2x^n(1 + x)\\
&=& x^2(x+1)^{2n+1} - 2x^{n+1} + (x^2+2x)^{n}(2x^2+3x).\quad\qed
\end{eqnarray*}

\begin{figure}[ht]
\hspace{4cm}
\includegraphics[width=7cm]{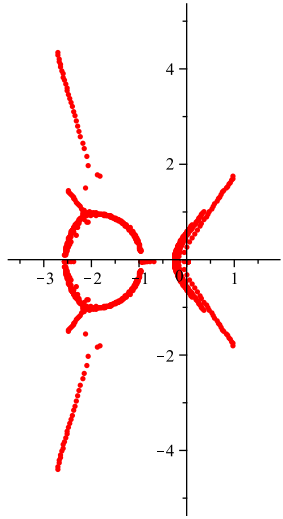}
\caption{\label{gbookroots} Domination roots of graphs $B_{n,5}$, for $1 \leq n \leq 30$.}
\end{figure}

\nt Figure~\ref{gbookroots} shows the domination roots of book graphs $B_{n,5}$ for $n \leq 30$.

\subsection{Limits of domination roots of book graphs $B_{n}$ and $B_{n,5}$.}

\nt In this section we consider  the complex domination roots of book graphs. The plot in Figures ~\ref{bookroots} and \ref{gbookroots} suggest that the roots tend to lie on a curve. In order to find the limiting curve, we will need a definition and a well known result.

\begin{newdefine}
If ${f_n(x)}$ is a family of (complex) polynomials, we say that a number $z \in \mathbb{C}$ is a limit of roots of ${f_n(x)}$ if either $f_n(z) = 0$ for all sufficiently large $n$ or z is a limit point of the set $\mathbb{R}({f_n(x)})$, where $\mathbb{R}({f_n(x)})$ is the union of the roots of the $f_n(x)$.
\end{newdefine}

\nt The following restatement of the Beraha-Kahane-Weiss theorem \cite{bkw}  can be found in \cite{brownhickman}.

\begin{teorem}\label{bkw}
Suppose ${f_n(x)}$ is a family of polynomials such that
\begin{eqnarray}
f_n(x) = \alpha_1(x)\lambda_1(x)^n + \alpha_2(x)\lambda_2(x)^n + ... + \alpha_k(x)\lambda_k(x)^n
\end{eqnarray}
where the $\alpha_i(x)$ and the $\lambda_i(x)$ are fixed non-zero polynomials, such that for no pair $i \neq j$ is $\lambda_i(x) \equiv \omega\lambda_j(x)$ for some $\omega \in \mathbb{C}$ of unit modulus. Then $z \in \mathbb{C}$ is a limit of roots of ${f_n(x)}$ if and only if either
\begin{itemize}
\item[(i)] two or more of the $\lambda_i(z)$ are of equal modulus, and strictly greater (in modulus) than the others; or
\item[(ii)] for some $j$, $\lambda_j(z)$ has modulus strictly greater than all the other $\lambda_i(z)$, and $\alpha_j(z) = 0.$
\end{itemize}
\end{teorem}

\nt The following Theorem gives the limits of the domination roots of book graphs $B_n$.

\begin{teorem}\label{mainroots}
The limit of domination roots of book  graphs are  $x= -\frac{1}{2}$ and $x=0$ together with
the part of the circle $|x+2|= 1$ with real part at least $\displaystyle{-\frac{3}{2}-\frac{\sqrt{2}}{2}}$, the portions of the hyperbola $(\Re (x) + 1)^{2} - (\Im (x))^{2} = \frac{1}{2}, ~\Re (x)\notin [\frac{-3-\sqrt{2}}{2},\frac{-2-\sqrt{2}}{2}]$, plus the portion of the curve $|x+1|^{2} = |x|$ with real part at most  $\displaystyle{-\frac{3}{2}-\frac{\sqrt{2}}{2}}$.
\end{teorem}

\nt{\bf Proof.} By Theorem~\ref{theorem6}, the domination polynomial of $B_n$ is,
\begin{eqnarray*}
D(B_n,x)&=&(2x+1)(x^2+2 x)^n + x^2(x+1)^{2n}- 2x^n\\
&=&\alpha_1(x)\lambda_1^n(x) + \alpha_2(x)\lambda_2^n(x) +\alpha_3(x)\lambda_3^n(x),
\end{eqnarray*}
\nt where
\[ \alpha_{1}(x) = 2x+1,~~\lambda_{1}(x) = x^2+2 x,\]
\[\alpha_{2}(x) =x^2,~~\lambda_{2}(x) = (x+1)^{2},\]
and
\[\alpha_{3}(x) = -2, ~~ \lambda_{3}(x) = x.\]

\nt Clearly $\alpha_1,~\alpha_2$ and $\alpha_3$ are not identically zero. Also, no $\lambda_{i} = \omega \lambda_{j}$
for $i\neq j$ and a  complex number $\omega $ of modulus $1$. Therefore, the initial conditions of Theorem~\ref{bkw} are satisfied. Now, applying part $(i)$ of Theorem~\ref{bkw}, we consider 4 different cases:

\begin{itemize}
\item[$(i)$] $|\lambda_1| = |\lambda_2| = |\lambda_3|$

\item[$(ii)$] $|\lambda_1| = |\lambda_2| > |\lambda_3|$

\item[$(iii)$] $|\lambda_1| = |\lambda_3| > |\lambda_2|$

\item[$(iv)$]  $|\lambda_2| = |\lambda_3| > |\lambda_1|$
\end{itemize}

\begin{figure}[ht]
\hspace{4cm}
\includegraphics[width=7cm]{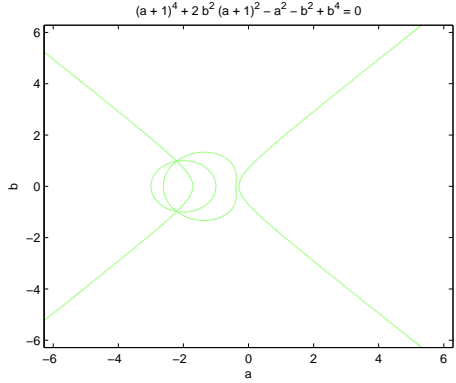}
\caption{\label{cardioid} The curves in  case $(i)$ in the proof of Theorem \label{mainroots}.}
\end{figure}

\nt {\bf case  $(i)$:} Assume that    $|x^2+2 x| = |(x+1)^{2}| = |x|$. Then   $|x^2+2 x| =  |x|$ implies that $x$ lies on the unit circle centered $-2~(|x - (-2)|= 1)$ and $|x^2+2 x| = |(x+1)^{2}|$ by setting $y=x+1$, that is,
\[ |y^{2}-1| = |y^{2}|.\]
To find this curve, let $a = \Re (y)$ and $b = \Im (y)$. Then by substituting in $y = a+ib$ and squaring both sides, we have
\[ (a^{2}-1-b^{2})^{2} + (2ab)^{2} =  (a^{2}-b^{2})^{2} + (2ab)^{2}.\]
This is equiavlent to
\[ a^{2} - b^{2} = \frac{1}{2},\]
a hyperbola. Hence, we converting back to variable $x$, we have the following hyperbola
\[(\Re (x) + 1)^{2} - (\Im (x))^{2} = \frac{1}{2}.\] Now suppose that
 $|(x+1)^{2}| = |x|$, this curve is semi-cardioid which has shown in Figure \ref{cardioid}.
 Therefore, the two points of intersection, $\frac{-3-\sqrt{2}}{2}\pm \frac{\sqrt{1+2\sqrt{2}}}{2}i$, are limits of roots.

\nt {\bf case  $(ii)$:} Assume that    $|x^2+2 x| = |(x+1)^{2}| > |x|$. Then $|x^2+2 x| = |(x+1)^{2}|$ implies that $x$ lies on the hyperbola $(\Re (x) + 1)^{2} - (\Im (x))^{2} = \frac{1}{2}$. And  $|x^2+2 x| > |x|$ implies that $x$ lies outside  the unit circle centered $-2~(|x - (-2)|= 1)$, and $|(x+1)^{2}| > |x|$ implies that $x$ lies outside  the curve $|(x+1)^{2}| = |x|$. Therefore, the complex numbers $x$ that satisfy
\[ (\Re (x) + 1)^{2} - (\Im (x))^{2} = \frac{1}{2}, \quad \Re (x)\notin [\frac{-3-\sqrt{2}}{2},\frac{-2-\sqrt{2}}{2}]\]
 are limits of roots.

\nt {\bf case  $(iii)$:} Assume that    $|x^2+2 x| =|x| > |(x+1)^{2}|$. Then  $|x^2+2 x| = |x|$ implies that $x$ lies on  the unit circle centered $-2~(|x - (-2)|= 1)$ and $|x^2+2 x| > |(x+1)^{2}|$ implies that $x$ satisfy in the following inequality
\[(\Re (x) + 1)^{2} - (\Im (x))^{2}< \frac{1}{2}.\]
The inequality  $|x| > |(x+1)^{2}|$  implies that $x$ lies inside  the curve $|(x+1)^{2}| = |x|$. Therefore, the complex numbers $x$ that satisfy  $|x - (-2)|= 1$ with real part at least $\frac{-3-\sqrt{2}}{2}$ are limits of roots.

\nt {\bf case  $(iv)$:} Assume that    $ |(x+1)^{2}| =|x| >|x^2+2 x|$. As we observed before,  the equality  $|(x+1)^{2}| = |x|$ is semi-cardioid which has shown in Figure \ref{cardioid}. The inequality  $|x| > |x^2+2 x|$ implies that $x$ lies inside the unit circle centered $-2~(|x - (-2)|= 1)$, and  $ |(x+1)^{2}| >|x^2+2 x|$ implies that $x$ satisfy in the following inequality
\[ (\Re (x) + 1)^{2} - (\Im (x))^{2} > \frac{1}{2}.\]
Therefore, the complex numbers $x$ that satisfy  on the curve $|(x+1)^{2}| = |x|$ with real part at most $\frac{-3-\sqrt{2}}{2}$ are limits of roots.

\nt  Finally by Part $(ii)$ of Theorem~\ref{bkw}, since  $\alpha_{3}$ is never $0$, and $\alpha_{2} = 0$ iff $x= 0$,  in this case $|\lambda_{2}(0)| = |1| > 0 = |\lambda_{1}(0)| = |\lambda_{3}(0)| $, and $\alpha_{1} = 0$ iff $x= -\frac{1}{2}$, and also in this case $|\lambda_{1}(-\frac{1}{2})| = |-\frac{3}{4}| > \frac{1}{4} = |\lambda_{2}(-\frac{1}{2})| $ and $|\lambda_{1}(-\frac{1}{2})| = |-\frac{3}{4}| > \frac{1}{2} = |\lambda_{3}(-\frac{1}{2})| $, so we conclude  $x= 0$ and  $x =-\frac{1}{2}$ are limit of domination roots of book graphs.

\nt The union of the curves and points above yield the desired result.\qquad\qed

\nt Along the same lines, we can  show:

\begin{teorem}
The limit of  roots of the domination polynomial of the book  graphs $B_{n,5}$, consist of
the part of the circle $|x+2|= 1$ with real part at least $\displaystyle{-\frac{3}{2}-\frac{\sqrt{2}}{2}}$, the portions of the hyperbola $(\Re (x) + 1)^{2} - (\Im (x))^{2} = \frac{1}{2}, ~\Re (x)\notin [\frac{-3-\sqrt{2}}{2},\frac{-2-\sqrt{2}}{2}]$, plus the portion of the curve $|x+1|^{2} = |x|$ with real part at most  $\displaystyle{-\frac{3}{2}-\frac{\sqrt{2}}{2}}$.
\end{teorem}


\bigskip


\begin{thebibliography}{99}


\bibitem{euro} S. Akbari, S. Alikhani and Y.H.  Peng, {\it Characterization of
graphs using domination polynomial}, Europ. J. Combin.,  Vol 31 (2010) 1714-1724.



\bibitem {saeid2}  S. Alikhani, {\it On the domination polynomials of non $P_4$-free graphs},
Iran. J. Math. Sci. Informatics, Vol. 8, No. 2 (2013) 49-55.


\bibitem{saeid1}  S. Alikhani, Y.H. Peng, {\it Introduction to domination polynomial of a graph}, Ars Combin., Vol. 114 (2014) pp. 257–266.

\bibitem{jason}  S. Alikhani, J.I. Brown, S. Jahari,  {\it On the domination polynomials of friendship graphs}, FILOMAT, to appear. Available at \texttt{http://arxiv.org/abs/1401.2092}.


\bibitem{bkw}
S.\ Beraha, J.\ Kahane, and N.\ Weiss,
{\it Limits of zeros of recursively defined families of polynomials}, in: G.\ Rota (Ed.),  Studies in foundations and combinatorics, Academic Press, New York, 1978, 213-232.

\bibitem{brownhickman}
J.I.\ Brown and C.A.\ Hickman,
{\it On chromatic roots of large subdivisions of graphs},
 Discr. Math. 242 (2002) 17-30.

\bibitem{brown} J.I. Brown and J. Tufts, {\it On the Roots of Domination Polynomials}, Graphs  Combin. 30 (2014), 527-547. doi: \texttt{10.1007/s00373-013-1306-z.}

\bibitem{Fruc} R. Frucht and F. Harary, {\it On the corona of two graphs},
Aequationes Mathematicae, vol. 4 (1970) 322-325.

\bibitem{garey}
M. R.\ Garey and D. S.\ Johnson, {\it Computers and Intractability:
A Guide to the Theory of $NP$-Completeness}, W. H. Freeman and Company,
New York, 1979



\bibitem{domination} T.W. Haynes, S.T. Hedetniemi, P.J. Slater, {\it Fundamentals of domination in graphs}, Marcel Dekker, NewYork, 1998.

\bibitem {Kot} T. Kotek, J. Preen, F. Simon, P. Tittmann, M. Trinks, {\it Recurrence relations and splitting formulas for the domination polynomial},  Elec. J. Combin. 19(3) (2012), \# P47.

\bibitem {Kotek} T. Kotek, J. Preen, P. Tittmann, {\it Domination polynomials of graph products}, arXive:1305.1475v2.



\bibitem{schi} Z. Ryj$\acute{a}\check{c}$ek, I. Schiermeyer, {\it The flower conjecture in special classes of graphs}, Discuss. Math. Graph Theory 15 (1995) 179–184.


\bibitem{Wal} M. Walsh,  {\it The hub number of a graph}, Int. J. Math. Comput. Sci., 1 (2006) 117-124.

\bibitem{htt} \texttt{http://mathworld. wolfram. com/DutchWindmillGraph. html}

\end{thebibliography}
\end{document}